\documentstyle{article}
\def\msc{{\bf M.S.C. 2000}:\ }
\def\kwd{{\bf Key words}:\ }
\title{On the Unification of Gravitational and Electromagnetic Fields}
\author{Kostadin Tren\v cevski}
\date{}
\begin{document}
\maketitle

\begin{abstract}
In the paper [4] is presented a theory which unifies
the gravitation theory and the mechanical effects, which
is different from the Riemannian theories like GTR. Moreover it
is built in the style of the electomagnetic field theory.
This paper is a continuation of [4] such that
the complex variant of that theory yields to the required
unification of gravitation and electromagnetism.
While the gravitational field is described by a scalar potential
$\mu $, taking a complex value of $\mu $ we obtain the unification theory.
For example the electric field appears to be imaginary 3-vector
field of acceleration, the magnetic field
appears to be imaginary 3-vector of angular velocity and the
imaginary part of a complex mass is just electric charge of the particle.
\end{abstract}
\msc 53B50, 53C80.\\
\kwd gravitational field, electromagnetic field, scalar potential.

\section{Preliminaries}
The present paper is upgrading of the paper [4] and so in 
this section will be repeated the basic concepts of the 
theory given in [4], which will need us in section 2. We convenient 
to use that $x^{4}=ict$ for the time coordinate and orthogonal matrix 
$A$ with complex elements will mean that $AA^{t}=I$. 

It this theory we deal with orthogonal frames of vector fields instead
of using curvilinear coordinates. In order this consideration to have
a physical meaning we will use coordinates $x^{i}$, but then 
$x^{i}$ $(i=1,2,3,4)$ are not functions of $x'^{j}$ $(j=1,2,3,4)$ if 
both of them are coordinates of noninertial systems, although 
the transition matrix $[\partial x^{i}/\partial x'^{j}]$ exists. 

In each infinitesimally small neighborhood of any point, the space-time continuum
can be considered as flat and the laws of Special Theory of Relativity
(STR) can be used there. Specially, the speed of the light 
measured in gravitational field in our laboratory, has the same 
value as measured in space without gravitation. But if 
we measure it on distance, such that there is a potential between us and 
the light ray, then we will not measure the same universal constant. 

The space-time in this theory is everywhere homogeneous which is not 
case with the General Theory of Relativity 
(GTR) where $g_{11}=g_{22}=g_{33}\neq \pm g_{44}$. The 
gravitation is described by one gravitation potential $\mu $. For
weak gravitational field one can assume that 
$\mu =1 + {\gamma M\over rc^{2}}$ where $\gamma $ is the universal 
gravitational constant. If $S$ and $S'$ are two systems which mutually rest
with parallel coordinate axes and assume that the observers from 
those systems are at points of differential gravitational potentials, then
the transformation between those two systems is given by
$$x'^{i} = \mu \cdot x^{i}\qquad (1\le i\le 4),\eqno{(1.1)}$$
where $\mu $ is the gravitational potential between the observers. 
The metric in general case is given by 
$g_{ij}=\mu ^{2}\delta_{ij}$ if we measure 
on distance, i.e. if there is a gravitational potential between us and 
object of measurement. If there is not such difference in the gravitational
potentials (if we measure in a laboratory), then the metric tensor is 
$g_{ij}=\delta _{ij}$. If two gravitational fields are given by the 
gravitational potentials $\mu _1 $ and $\mu _2 $, then the superposition
of those two fields is a gravitational field with potential $\mu _1 \mu _2 $.
In special case if $\mu_2 =C$ is constant in the space-time, then all
physical laws should be the same as for the field given by $\mu_1 $. 

This theory is built on the basis of 4-vector potential $U_i $. 
If the magnitude of this vector changes then 
there exists a gravitational field and if the direction of this vector 
changes then there exist some mechanical forces, for example centrifugal. 
Indeed, this field can be written in general case in the form 
$$ (U_{1},U_{2},U_{3},U_{4}) = 
{\mu \over \sqrt {1-u^{2}/c^{2}}}\Bigl ( 
{u_{1}\over ic}, 
{u_{2}\over ic}, 
{u_{3}\over ic}, 
1\Bigr ) ,\eqno{(1.2)}$$
where 
${\bf u}=(u_{1},u_{2},u_{3})$ is 3-vector of velocity and $\mu $ is the 
gravitational potential. 

{\it Example 1.} The STR is characterized by $\mu =const.$ and 
${\bf u}$ is a constant vector field. 

{\it Example 2.} The gravitational field of a planet 
which does not rotate and does not move is given by the vector field 
$(0,0,0,\mu )$. 

{\it Example 3.} If we measure with respect to a noninertial system 
which rotates
around the $x^{3}$ axis with angular velocity $w$, then the 4-vector 
potential obtains by putting $\mu =1$ and 
${\bf u}=(x^{2}\omega ,-x^{1}\omega ,0)$. 

Note that the main problem is to determine the 4-vector potential.  
In case of gravitational field, ${\bf u}$ shows the 
3-vector of velocity of the mass which produce the gravitation. 
As we will see
later it is always very convenient if we choose a system 
such that ${\bf u}$ is zero vector at the considered point. 

Next we introduce a non-linear connection induced by the 
4-vector potential $U_i $. Let denote by $N^{i}_{j}ds$ the form of 
connection in the direction of the unit vector $V^{i}=dx^{i}/ds$. In order
the metric tensor $g_{ij}=\mu ^{2}\delta _{ij}$ to be parallel, i.e. 
$dg_{ij}/ds - g_{ir}N^{r}_{j}-g_{rj}N^{r}_{i}=0$, we put 
$N^{i}_{j}=\delta ^{i}_{j}d(ln \mu )/ds + S^{i}_{j}$ and obtain that 
$S^{i}_{j}=-S^{j}_{i}$. Note that $S^{i}_{j}ds$ is a form of connection 
which makes the metric tensor field $g_{ij}=\delta _{ij}$ parallel. For any
vector field $A^{i}$ one verifies the identity 
$$dA^{i}/ds + N^{i}_{j}A^{j} = {1\over \mu }[d(A^{i}\mu )/ds + S^{i}_{j}
(A^{j}\mu )]\eqno{(1.3)}$$
and $d(A^{i}\mu )/ds + S^{i}_{j}(A^{j}\mu )$ is the covariant derivative 
with respect to the antisymmetric connection $S^{i}_{j}$. The previous 
formula shows that the covariant derivatives with respect to
these two connections differ only by their magnitudes which should be 
covariant with respect to the corresponding metrics. This permits us 
to consider further only the connection $S^{i}_{j}$ and to assume that
$g_{ij}=\delta _{ij}$ at the considered point. Thus in future there will be
no difference between upper and lower indices. 

Now let us introduce the following tensor 
$$F_{ij}={1\over \mu }(\partial U_{i}/\partial x^{j}- 
\partial U_{j}/\partial x^{i}).\eqno{(1.4)}$$
If we denote by $V_{i}$ the 4-vector of velocity of the 
experimental particle, i.e. 
$$V=\Bigl ( 
{v_{1}\over ic\sqrt {1-v^{2}/c^{2}}}, 
{v_{2}\over ic\sqrt {1-v^{2}/c^{2}}}, 
{v_{3}\over ic\sqrt {1-v^{2}/c^{2}}}, 
{1\over \sqrt {1-v^{2}/c^{2}}}\Bigr ) \eqno{(1.5)}$$
one can easily verify that the following simple equation
$$ dV_{i}/ds = F_{ij}V_{j}\eqno{(1.6)}$$
which is analogous to the Lorentz force in the electromagnetic theory, 
contains all mechanical accelerations (Coriolis acceleration, 
centrifugal and transversely acceleration) and gravitation accelerations
for Newton approach. In order to obtain more precise theory, we consider 
the following two tensors 
$$\phi _{ij}={1\over 2}F_{ij} - {1\over 2}\mu ^{-2}
(U_{i}U_{k}F_{jk}- U_{j}U_{k}F_{ik}),\eqno{(1.7)}$$
$$P_{ij}=\delta _{ij} - {1\over \mu + U_{s}V_{s}}
(\mu V_{i}V_{j}+V_{i}U_{j}+U_{i}V_{j}+{1\over \mu }U_{i}U_{j}) + 
{2\over \mu }U_{j}V_{i}.\eqno{(1.8)}$$
These two tensors together with $F_{ij}$ are invariant under the gauge 
transformation $\mu \rightarrow C\mu $ where $C$ is a constant. The tensors 
$F_{ij}$ and $\phi _{ij}$ are antisymmetric and $P_{ij}$ is an orthogonal 
matrix. It can be verified by using the identities 
$U_{i}U_{i}=\mu ^{2}$ and $V_{i}V_{i}=1$. Now the required connection 
$S_{ij}$ is introduced by
$$S_{ij} = - P_{ri}\phi _{rk}P_{kj},\eqno{(1.9)}$$
or in matrix form $S=-P^{t}\phi P$. Although this connection seems 
to be very complicated, now we will show the opposite. Indeed, assume that 
$(U_{i})=(0,0,0,1)$ at the considered point. Then 
$$ F = \left [\matrix{
0 & -2i\omega _{3}/c & 2i\omega _{2}/c & -a_{1}/c^{2}\cr 
2i\omega _{3}/c & 0 & -2i\omega _{1}/c & -a_{2}/c^{2}\cr 
-2i\omega _{2}/c & 2i\omega _{1}/c & 0 & -a_{3}/c^{2}\cr 
a_{1}/c^{2} & a_{2}/c^{2} & a_{3}/c^{2} & 0\cr }\right ], \eqno{(1.10)}$$
$$ \phi = \left [\matrix{
0 & -i\omega _{3}/c & i\omega _{2}/c & -a_{1}/c^{2}\cr 
i\omega _{3}/c & 0 & -i\omega _{1}/c & -a_{2}/c^{2}\cr 
-i\omega _{2}/c & i\omega _{1}/c & 0 & -a_{3}/c^{2}\cr 
a_{1}/c^{2} & a_{2}/c^{2} & a_{3}/c^{2} & 0\cr }\right ], \eqno{(1.11)}$$
where 
$$ {\bf a}=(a_{1},a_{2},a_{3})=c^{2}grad U_{4} + 
\partial {\bf u}/\partial t \quad \hbox { and }\quad
{\bf w}= -{1\over 2}rot {\bf u}\eqno{(1.12)}$$
represent the 3-vector of acceleration and the 3-vector of angular velocity. 
In that special case $(U_{i})=(0,0,0,1)$ the tensor $P$ is given by 
$$ P = \left [\matrix{
1-{1\over \nu }V_{1}^{2} & -{1\over \nu }V_{1}V_{2} &
-{1\over \nu }V_{1}V_{3} & V_{1}\cr 
-{1\over \nu }V_{2}V_{1} & 1-{1\over \nu }V_{2}^{2} &
-{1\over \nu }V_{2}V_{3} & V_{2}\cr 
-{1\over \nu }V_{3}V_{1} & -{1\over \nu }V_{3}V_{2} & 
1-{1\over \nu }V_{3}^{2} & V_{3}\cr 
-V_{1} & -V_{2} & -V_{3} & V_{4}\cr }\right ], \eqno{(1.13)}$$
where $V_{1}, V_{2}, V_{3}, V_{4}$ were given previously and 
$\nu =1+V_{4}$. Note that the matrix $P$ represents just a Lorentz 
transformation for the 3-vector $-{\bf v}=(-v_{1},-v_{2},-v_{3})$. Now the 
connection $S_{ij}$ becomes much more clear and acceptable. 

In [4] are given some consequences of the connection $S_{ij}$. Here we note
that for the angle of deflection of light ray it is obtained 
$\Delta \alpha ={4\gamma M\over Rc^{2}}$ which is the same as in the GTR
and also for the angle between two perihelions of a planet it is obtained 
the same value as in the GTR. In [4] is also discussed the paradox of the 
twins in gravitational field and noninertial systems. 

Note that the introduced connection $S_{ij}$ is non-linear. The first step 
of approximation yields to the following linear connection 
$\Gamma _{k}=(\Gamma ^{i}_{jk})$ such that $S\approx \Gamma _{s}V_{s}$: 
$$ \Gamma _{1} = \left [ \matrix{
0 & a_{2}/c^{2} & a_{3}/c^{2} & 0\cr
-a_{2}/c^{2} & 0 & 0 & -iw_{3}/c \cr 
-a_{3}/c^{2} & 0 & 0 & iw_{2}/c \cr 
0 & iw_{3}/c & -iw_{2}/c &0\cr }\right ],\eqno{(1.14a)} $$
$$ \Gamma _{2} = \left [ \matrix{
0 & -a_{1}/c^{2} & 0 &  iw_{3}/c\cr
a_{1}/c^{2} & 0 & a_{3}/c^{2} & 0\cr
0 & -a_{3}/c^{2} & 0 & -iw_{1}/c \cr 
-iw_{3}/c & 0 & iw_{1}/c & 0\cr }\right ],\eqno{(1.14b)} $$
$$ \Gamma _{3} = \left [ \matrix{
0 & 0 & -a_{1}/c^{2} & -iw_{2}/c\cr
0 & 0 & -a_{2}/c^{2} & iw_{1}/c\cr
a_{1}/c^{2} & a_{2}/c^{2} & 0 & 0 \cr 
iw_{2}/c & -iw_{1}/c &0&0\cr }\right ], \eqno{(1.14c)}$$
$$ \Gamma _{4} = \left [ \matrix{
0 & iw_{3}/c & -iw_{2}/c & a_{1}/c^{2}\cr 
-iw_{3}/c & 0 & iw_{1}/c & a_{2}/c^{2}\cr 
iw_{2}/c & -iw_{1}/c & 0 & a_{3}/c^{2}\cr 
-a_{1}/c^{2} & -a_{2}/c^{2} & -a_{3}/c^{2} & 0 \cr }\right ],\eqno{(1.14d)}$$
where we assumed that $(U_{i})=(0,0,0,1)$ at the considered point. 
Note that this approximation is of order like  the  GTR  and  the Newton's 
approximation obtains by replacing with zero the components 
$\Gamma ^{i}_{jk}$ for $i,j,k\in \{1,2,3\}$. Using this linear connection 
the curvature tensor is easy to be calculated and in [4] it is shown that
the Einstein equations of gravitational field 
$$ R_{ij} - {1\over 2} g_{ij}R = 8\gamma \pi c^{-2} T_{ij}\eqno{(1.15)}$$
are satisfied, where the energy-momentum tensor $T_{ij}$ is given by 
$T_{ij}=\rho U_{i}U_{j}$. It is also shown in [4] that the following 
equations 
$$(\mu F_{ij})_{;k}+(\mu F_{jk})_{;i}+(\mu F_{ki})_{;j}=0 ,\eqno{(1.16)}$$
$$ (\mu F_{ij})_{;j} = 4\gamma \pi c^{-2} J_{i} ,\eqno{(1.17)}$$
analogous to the Maxwell equations, also hold at that order of 
approximation. Here $J_{i}=\rho U_{i}$
is the 4-current density for gravitation, analogous to the 4-current 
density for the electromagnetism. 

\section{Unification of gravitation and electromagnetism}

Note that the tensor $F_{ij}$ defined by (1.4) in special case 
when $U=(0,0,0,1)$ at the considered point is given by $(1.10)$ where 
${\bf a}$ and ${\bf w}$ are the 3-vectors of acceleration and angular
velocity. The Newton approximation of motion of the particle in such 
field is given by (1.6) which is analogous to the Lorentz force in 
electromagnetic field. Thus the tensor $F_{ij}$ has the same meaning 
as the tensor of electromagnetic field and the vector field $U_{i}$ has 
the same meaning like the 4-vector potential $A_{i}$ for electromagnetic 
field. This permits us that gravitation theory is analogous to 
the electromagnetic theory and to unify them as follows. 

The previous discussion associates us to consider the electric field 
as imaginary acceleration, and the  
magnetic field is imaginary angular velocity, i.e. 
$${\bf E} \rightarrow i\lambda ^{-1}{\bf a}\qquad \hbox { and }\qquad  
{\bf H} \rightarrow i\lambda ^{-1}2c{\bf w},\eqno{(2.1)}$$
where $\lambda $ is an universal constant. Comparing the formulas for 
forces ${\bf f}=m{\bf a}$ (if ${\bf w}=0$) and 
${\bf f}=e{\bf E}$ (if ${\bf H}=0$), 
we notice that the electricity $e$
is an imaginary mass, i.e. $e\rightarrow -i\lambda m$ or 
$e\rightarrow i\lambda m$, which depends on the choice of the charge 
$e$: positive or negative. In order to 
determine the universal constant $\lambda $, we compare the Newton's force
$f=\gamma m_{1}m_{2}/r^{2}$
and the Coulon's force $f={1\over 4\pi \epsilon _{0}}e_{1}e_{2}/r^{2}$
and obtain that 
$\lambda =\sqrt{4\pi \epsilon _{0} \gamma }$. 
Hence, if we have a particle with
electricity $e$ and mass $m$, it has complex mass $M=m-i{1\over \lambda }e$. 

Note that the gravitational potential $\mu =1+{\gamma M\over rc^{2}}$ 
has complex value if M has the previous complex mass. Indeed, we accept that
in absence of the gravitational field, in a neighborhood of particle with 
electric charge $e$, the potential is 
$$ \nu = 1 - {ie\over \lambda rc^{2}}.\eqno{(2.2)}$$
Note that according to the law of superposition of potentials, if
$\nu _{1}=1+iV_{1}$, $\nu _{2}=1+iV_{2}$ are two electromagnetic potentials, 
then 
$$\nu _{1}\nu _{2}=1-V_{1}V_{2} + i(V_{1}+V_{2}) = (1-V_{1}V_{2})
\Bigl (1+ i{V_{1}+V_{2}\over 1-V_{1}V_{2}}\Bigr )$$
and hence the superposition gives an electromagnetic potential 
$\nu =1+iV$ where $V = {V_{1}+V_{2}\over 1-V_{1}V_{2}}$ and also a 
gravitational potential $\mu =1-V_{1}V_{2}$, which may be neglected if 
$V_{1}\approx 0$ and $V_{2}\approx 0$. 

In the unified theory the 4-vector potential is 
a vector field including the gravitational and electromagnetic potentials 
and it is given by 
$$ (U_{1},U_{2},U_{3},U_{4}) = 
{\mu \nu \over \sqrt {1-u^{2}/c^{2}}}\Bigl ( 
{u_{1}\over ic}, 
{u_{2}\over ic}, 
{u_{3}\over ic},
1\Bigr ) ,\eqno{(2.3)}$$
where $\mu \in R$ and ${\nu -1\over i}\in R$. Note that 
$${u_{1}\over c} = {iU_{1}\over U_{4}}\in R, \quad 
{u_{2}\over c} = {iU_{2}\over U_{4}}\in R, \quad 
{u_{3}\over c} = {iU_{3}\over U_{4}}\in R. $$
Conversely, if $(U_{1},U_{2},U_{3},U_{4})$ is any vector field such that 
${iU_{1}\over U_{4}}\in R$, ${iU_{2}\over U_{4}}\in R$, 
${iU_{3}\over U_{4}}\in R$, then $u_{1}, u_{2}, u_{3}, \mu ,\nu $ are 
uniquely determined. 
\par 
Further, according to the first step of approximation the classical results 
will be obtained. Indeed, let us assume that $\mu \approx 1$, 
$\nu \approx 1$. Indeed, we use that $\mu =\nu =1$, but 
${\partial \mu \over \partial x^{j}}\neq 0$ and 
${\partial \nu \over \partial x^{j}}\neq 0$. Analogously to (1.4) we 
obtain the antisymmetric tensor 
$$ \Phi _{pq} = {\partial U_{p}\over \partial x^{q}} -
{\partial U_{q}\over \partial x^{p}}.\eqno{(2.4)}$$ 
By decomposing it into real and imaginary part we obtain 
$$ \Phi _{pq} = F_{pq} + i\lambda \Psi _{pq} ,\eqno{(2.5)}$$
where $\Psi _{pq}$ is tensor proportional with the tensor of 
electromagnetic field and it is given by 
$$  \Psi = \left [ \matrix { 
0 & iH_{3}/c^{2} & -iH_{2}/c^{2} & E_{1}/c^{2}\cr   
-iH_{3}/c^{2} & 0 & iH_{1}/c^{2} & E_{2}/c^{2}\cr 
iH_{2}/c^{2} & -iH_{1}/c^{2} & 0 & E_{3}/c^{2}\cr 
-E_{1}/c^{2}&-E_{2}/c^{2}&-E_{3}/c^{2}&0\cr }
\right ].  \eqno{(2.6)} $$
Note that the tensor $\Psi $ can be written as 
$$ \Psi _{pq} = {i\over c^{2}} \Bigl ({\partial A_{q}\over \partial x^{p}}
-{\partial A_{p}\over \partial x^{q}}\Bigr ) ,\eqno{(2.7)}$$
where 
$$ A_{1} = i{u_{1}c\mu (Im\nu )\over \sqrt {1-u^{2}/c^{2}}\lambda }, \quad 
A_{2} = i{u_{2}c\mu (Im\nu )\over \sqrt {1-u^{2}/c^{2}}\lambda }, \quad
A_{3} = i{u_{3}c\mu (Im\nu )\over \sqrt {1-u^{2}/c^{2}}\lambda }, $$
$$A_{4} = -{c^{2}\mu (Im\nu )\over \sqrt {1-u^{2}/c^{2}}\lambda }, 
\eqno{(2.8)}$$
which gives the well known 4-vector of potential of the 
electromagnetic field. 

Now let us consider the motion of a charged particle in the field 
$\Phi _{ij}$. Let the particle has mass $m$ and electricity $e$, i.e. 
its complex mass is $M=m-{i\over \lambda }e$. Then analogously to the 
equation of motion (1.6) and at the some level of approximation it holds 
$$ F_{p} = M \Phi _{pq}V_{q} ,\eqno{(2.9)}$$
where $V_{q}$ is 4-vector velocity of the particle and $F_{p}$
is 4-vector of force. If we exchange the value of M, the values of 
$\Phi _{pq}$, the real part of the last equation of motion contains 
both action of the mechanical force and gravitational and also the force of 
charged particle in electromagnetic field and it is the force which we 
measure by experiments and which determines the motion of the particle.

Finally note that the impulse of the considered particle can be 
calculated by 
$$Re[(m-i{1\over \lambda }e)({\bf v} - {i\over c}\lambda {\bf A})] 
= m{\bf v} - e{\bf A}/c ,
\eqno{(2.10)}$$
where ${\bf v}$ and ${\bf A}$ are the 3-vectors of velocity and 
electromagnetic potential and it is in accordance 
with the accepted expression for impulse. 

For the sake of simplicity we showed the unification on the Newton's 
level of approximation. For more precise trajectory of motion of 
particle one should use the connection (1.9) or at least
the linear connection (1.14) in complex case.

\medskip

\noindent Institute of Mathematics 
\par 
\noindent St. Cyril and Methodius University
\par 
\noindent P.O.Box 162, 1000 Skopje, Macedonia 
\par 
\noindent e-mail: kostatre@iunona.pmf.ukim.edu.mk

\begin{thebibliography}{9}
\bibitem{L}{\sc D.F.Lawden}, Tensor Calculus and Relativity, Chapmann
	and Hall, 1976.
\bibitem{M}{\sc K.B.Marathe, M.Modugno}, Polynomial Connections on Affine
	Bundles, Tensor N.S., 50 (1991), 35-46.
\bibitem{Mi}{\sc C.W.Misner, K.S.Thorne, J.A.Wheeler}, Gravitation, W.H.
	Freeman and Co., San Francisco, 1973.
\bibitem{paper} {\sc K.Tren\v cevski}, One Model of Gravitation  and
	Mechanics, Tensor N.S., 53 (1993), 70-82.
\end{thebibliography}
\end{document}